\documentclass[11pt,english]{article}
\usepackage[T1]{fontenc}
\usepackage[latin9]{inputenc}
\usepackage{geometry}
\usepackage{prettyref}
\usepackage{amsthm}
\usepackage{amsmath}
\usepackage{amssymb}
\usepackage{setspace}
\usepackage{esint}
\makeatletter
\theoremstyle{plain}
\newtheorem{thm}{\protect\theoremname}
  \theoremstyle{definition}
  \newtheorem{example}[thm]{\protect\examplename}
  \theoremstyle{plain}
  \newtheorem{lem}[thm]{\protect\lemmaname}
  \theoremstyle{plain}
  \newtheorem{prop}[thm]{\protect\propositionname}

\newcounter{qcounter}

\newrefformat{pro}{Proposition \ref{#1}}

\makeatother

\usepackage{babel}
  \providecommand{\examplename}{Example}
  \providecommand{\lemmaname}{Lemma}
  \providecommand{\propositionname}{Proposition}
\providecommand{\theoremname}{Theorem}

\begin{document}

\title{On $\varphi$-families of probability distributions}

\author{Rui F.\ Vigelis%
\thanks{Computer Engineering, Campus Sobral, Federal University of Ceará,
Sobral-CE, Brazil. E-mail: rfvigelis@ufc.br.%
},\and Charles C.\ Cavalcante%
\thanks{Wireless Telecommunication Research Group, Department of Teleinformatics
Engineering, Federal University of Ceará, Fortaleza-CE, Brazil. E-mail:
charles@ufc.br.%
}}
\maketitle
\begin{abstract}
We generalize the exponential family of probability distributions.
In our approach, the exponential function is replaced by a $\varphi$-function,
resulting in a $\varphi$-family of probability distributions. We
show how $\varphi$-families are constructed. In a $\varphi$-family,
the analogue of the cumulant-generating function is a normalizing
function. We define the $\varphi$-divergence as the Bregman divergence
associated to the normalizing function, providing a generalization
of the Kullback--Leibler divergence. A formula for the $\varphi$-divergence
where the $\varphi$-function is the Kaniadakis' $\kappa$-exponential
function is derived.
\end{abstract}

\section{Introduction}

Let $(T,\Sigma,\mu)$ be a $\sigma$-finite, non-atomic measure space.
We denote by $\mathcal{P}_{\mu}=\mathcal{P}(T,\Sigma,\mu)$ the family
of all probability measures on $T$ that are equivalent to the measure
$\mu$. The probability family $\mathcal{P}_{\mu}$ can be represented
as (we adopt the same symbol $\mathcal{P}_{\mu}$ for this representation)
\[
\mathcal{P}_{\mu}=\{p\in L^{0}:p>0\text{ and }\mathbb{E}[p]=1\},
\]
where $L^{0}$ is the linear space of all real-valued, measurable
functions on $T$, with equality $\mu$-a.e., and $\mathbb{E}[\cdot]$
denotes the expectation with respect to the measure $\mu$. 

The family $\mathcal{P}_{\mu}$ can be equipped with a structure of
$C^{\infty}$-Banach manifold, using the Orlicz space $L^{\Phi_{1}}(p)=L^{\Phi_{1}}(T,\Sigma,p\cdot\mu)$
associated to the Orlicz function $\Phi_{1}(u)=\exp(u)-1$, for $u\geq0$.
With this structure, $\mathcal{P}_{\mu}$ is called the \textit{exponential
statistical manifold}, whose construction was proposed in \cite{Pistone:1995}
and developed in \cite{Pistone:1999,Cena:2007,Grasselli:2010}. Each
connected component of the exponential statistical manifold gives
rise to an \textit{exponential family of probability distributions
$\mathcal{E}_{p}$} (for each $p\in\mathcal{P}_{\mu}$). Each element
of $\mathcal{E}_{p}$ can be expressed as 
\begin{equation}
\boldsymbol{e}_{p}(u)=e^{u-K_{p}(u)}p,\qquad\text{for }u\in\mathcal{B}_{p},\label{eq:e_p}
\end{equation}
for a subset $\mathcal{B}_{p}$ of the Orlicz space $L^{\Phi_{1}}(p)$.
$K_{p}$ is the cumulant-generating functional $K_{p}(u)=\log\mathbb{E}_{p}[e^{u}]$,
where $\mathbb{E}_{p}[\cdot]$ is the expectation with respect to
$p\cdot\mu$. If $c$ is a measurable function such that $p=e^{c}$,
then \prettyref{eq:e_p} can be rewritten as 
\begin{equation}
\boldsymbol{e}_{p}(u)=e^{c+u-K_{p}(u)\cdot\boldsymbol{1}_{T}},\qquad\text{for }u\in\mathcal{B}_{p},\label{eq:e_pc}
\end{equation}
where $\boldsymbol{1}_{A}$ is the indicator function of a subset
$A\subseteq T$. A generalization of expression \prettyref{eq:e_p}
was given in \cite{Pistone:2009}, where the exponential function
is replaced by a $\kappa$-exponential function. In our generalization,
we make use of expression \prettyref{eq:e_pc}.

In the \textit{$\varphi$-family of probability distributions $\mathcal{F}_{c}^{\varphi}$},
which we propose, the exponential function is replaced by the so called
\textit{$\varphi$-function} $\varphi\colon T\times\overline{\mathbb{R}}\rightarrow[0,\infty]$.
The function $\varphi(t,\cdot)$ has a ``shape'' which is similar
to that of an exponential function, with an arbitrary rate of increasing.
For example, we found that the $\kappa$-exponential function satisfies
the definition of $\varphi$-functions. As in the exponential family,
the $\varphi$-families are the connected component of $\mathcal{P}_{\mu}$,
which is endowed with a structure of $C^{\infty}$-Banach manifold,
using $\varphi$ in the place of an exponential function. Let $c$
be any measurable function such that $\varphi(t,c(t))$ belongs to
$\mathcal{P}_{\mu}$. The elements of the $\varphi$-family of probability
distributions $\mathcal{F}_{c}^{\varphi}$ are given by 
\begin{equation}
\boldsymbol{\varphi}_{c}(u)(t)=\varphi(t,c(t)+u(t)-\psi(u)u_{0}(t)),\qquad\text{for }u\in\mathcal{B}_{c}^{\varphi},\label{eq:phi_c}
\end{equation}
for a subset $\mathcal{B}_{c}^{\varphi}$ of a Musielak--Orlicz space
$L_{c}^{\varphi}$. The \textit{normalizing function} $\psi\colon\mathcal{B}_{c}^{\varphi}\rightarrow[0,\infty)$
and the measurable function $u_{0}\colon T\rightarrow[0,\infty)$
in \prettyref{eq:phi_c} replaces $K_{p}$ and $\boldsymbol{1}_{T}$
in \prettyref{eq:e_pc}, receptively. The function $u_{0}$ is not
arbitrary. In the text, we will show how $u_{0}$ can be chosen.

We define the \textit{$\varphi$-divergence} as the a Bregman divergence
associated to the normalizing function $\psi$, providing a generalization
of the Kullback--Leibler divergence. Then geometrical aspects related
to the $\varphi$-family can be developed, since the Fisher information
(from which the Information Geometry \cite{Amari:2000,Murray:1993}
is based) is derived from the divergence. A formula for the $\varphi$-divergence
where the $\varphi$-function is the Kaniadakis' $\kappa$-exponential
function \cite{Kaniadakis:2002,Naudts:2011} is derived, which we
called the \textit{$\kappa$-divergence}.

We expect that an extension of our work will provide advances in other
areas, like in Information Geometry or in the non-parametric, non-commutative
setting \cite{Gibilisco:2010,Petz:2008}. The rest of this paper is
organized as follows. \prettyref{sec:musielak_orlicz_spaces} deals
with the topics of Musielak--Orlicz spaces we will use in the the
construction of the $\varphi$-family of probability distributions.
In \prettyref{sec:exponential_statistical_manifold}, the exponential
statistical manifold is reviewed. The construction of the $\varphi$-family
of probability distributions is given in \prettyref{sec:construction_phi_family}.
Finally, the $\varphi$-divergence is derived in \prettyref{sec:divergence}.

\section{Musielak--Orlicz spaces\label{sec:musielak_orlicz_spaces}}

In this section we provide a brief introduction to Musielak--Orlicz
(function) spaces, which are used in the construction of the exponential
and $\varphi$-families. A more detailed exposition about these spaces
can be found in \cite{Musielak:1983,Krasnoselskii:1961,Rao:1991}.

We say that $\Phi\colon T\times[0,\infty]\rightarrow[0,\infty]$ is
a \textit{Musielak--Orlicz function} when, for $\mu$-a.e.\ $t\in T$,
\begin{itemize}
\item [(i)] $\Phi(t,\cdot)$ is convex and lower semi-continuous,
\item [(ii)] $\Phi(t,0)=\lim_{u\downarrow0}\Phi(t,u)=0$ and $\Phi(t,\infty)=\infty$,
\item [(iii)] $\Phi(\cdot,u)$ is measurable for all $u\geq0$.
\end{itemize}
Items (i)--(ii) guarantee that $\Phi(t,\cdot)$ is not equal to $0$
or $\infty$ on the interval $(0,\infty)$. A Musielak--Orlicz function
$\Phi$ is said to be an \textit{Orlicz function} if the functions
$\Phi(t,\cdot)$ are identical for $\mu$-a.e.\ $t\in T$.

Define the functional $I_{\Phi}(u)=\int_{T}\Phi(t,|u(t)|)d\mu$, for
any $u\in L^{0}$. The \textit{Musielak--Orlicz space}, \textit{Musielak--Orlicz
class}, and \textit{Morse--Transue space}, are given by 
\begin{align*}
L^{\Phi} & =\{u\in L^{0}:I_{\Phi}(\lambda u)<\infty\text{ for some }\lambda>0\},\\
\tilde{L}^{\Phi} & =\{u\in L^{0}:I_{\Phi}(u)<\infty\},\\
\intertext{and}E^{\Phi} & =\{u\in L^{0}:I_{\Phi}(\lambda u)<\infty\text{ for all }\lambda>0\},
\end{align*}
respectively. If the underlying measure space $(T,\Sigma,\mu)$ have
to be specified, we write $L^{\Phi}(T,\Sigma,\mu)$, $\tilde{L}^{\Phi}(T,\Sigma,\mu)$
and $E^{\Phi}(T,\Sigma,\mu)$ in the place of $L^{\Phi}$, $\tilde{L}^{\Phi}$
and $E^{\Phi}$, respectively. Clearly, $E^{\Phi}\subseteq\tilde{L}^{\Phi}\subseteq L^{\Phi}$.
The Musielak--Orlicz space $L^{\Phi}$ can be interpreted as the smallest
vector subspace of $L^{0}$ that contains $\tilde{L}^{\Phi}$, and
$E^{\Phi}$ is the largest vector subspace of $L^{0}$ that is contained
in $\tilde{L}^{\Phi}$.

The Musielak--Orlicz space $L^{\Phi}$ is a Banach space when it is
endowed with the \textit{Luxemburg norm} 
\[
\Vert u\Vert_{\Phi}=\inf\Bigl\{\lambda>0:I_{\Phi}\Bigl(\frac{u}{\lambda}\Bigr)\leq1\Bigr\},
\]
or the \textit{Orlicz norm}
\[
\Vert u\Vert_{\Phi,0}=\sup\biggl\{\biggl|\int_{T}uvd\mu\biggr|:v\in\tilde{L}^{\Phi^{*}}\text{ and }I_{\Phi^{*}}(v)\leq1\biggr\},
\]
where $\Phi^{*}(t,v)=\sup\nolimits _{u\geq0}(uv-\Phi(t,u))$ is the
\textit{Fenchel conjugate} of $\Phi(t,\cdot)$. These norms are equivalent
and the inequalities $\Vert u\Vert_{\Phi}\leq\Vert u\Vert_{\Phi,0}\leq2\Vert u\Vert_{\Phi}$
hold for all $u\in L^{\Phi}$.

If we can find a non-negative function $f\in\tilde{L}^{\Phi}$ and
a constant $K>0$ such that
\[
\Phi(t,2u)\leq K\Phi(t,u),\quad\text{for all }u\geq f(t),
\]
then we say that $\Phi$ satisfies the \textit{$\Delta_{2}$-condition},
or belong to the \textit{$\Delta_{2}$-class} (denoted by $\Phi\in\Delta_{2}$).
When the Musielak--Orlicz function $\Phi$ satisfies the $\Delta_{2}$-condition,
$E^{\Phi}$ coincides with $L^{\Phi}$. On the other hand, if $\Phi$
is finite-valued and does not satisfy the $\Delta_{2}$-condition,
then the Musielak--Orlicz class $\tilde{L}^{\Phi}$ is not open and
its interior coincides with 
\[
B_{0}(E^{\Phi},1)=\{u\in L^{\Phi}:\inf_{v\in E^{\Phi}}\Vert u-v\Vert_{\Phi,0}<1\},
\]
or, equivalently, $B_{0}(E^{\Phi},1)\varsubsetneq\tilde{L}^{\Phi}\varsubsetneq\overline{B}_{0}(E^{\Phi},1)$.

\section{The exponential statistical manifold\label{sec:exponential_statistical_manifold}}

This section starts with the definition of a $C^{k}$-Banach manifold
\cite{Lang:1995}. A \textit{$C^{k}$-Banach manifold} is a set $M$
and a collection of pairs $(U_{\alpha},\boldsymbol{x}_{\alpha})$
($\alpha$ belonging to some indexing set), composed by open subsets
$U_{\alpha}$ of some Banach space $X_{\alpha}$, and injective mappings
$\boldsymbol{x}_{\alpha}\colon U_{\alpha}\rightarrow M$, satisfying
the following conditions: \begin{list}{(bm\arabic{qcounter})}{\usecounter{qcounter}\setlength{\leftmargin}{3.5em}\setlength{\labelwidth}{3.5em}}
\item the sets $\boldsymbol{x}_{\alpha}(U_{\alpha})$ cover $M$, i.e.,
$\bigcup_{\alpha}\boldsymbol{x}_{\alpha}(U_{\alpha})=M$; 
\item for any pair of indices $\alpha,\beta$ such that $\boldsymbol{x}_{\alpha}(U_{\alpha})\cap\boldsymbol{x}_{\beta}(U_{\beta})=W\neq\emptyset$,
the sets $\boldsymbol{x}_{\alpha}^{-1}(W)$ and $\boldsymbol{x}_{\beta}^{-1}(W)$
are open in $X_{\alpha}$ and $X_{\beta}$, respectively; and 
\item the \textit{transition map} $\boldsymbol{x}_{\beta}^{-1}\circ\boldsymbol{x}_{\alpha}\colon\boldsymbol{x}_{\alpha}^{-1}(W)\rightarrow\boldsymbol{x}_{\beta}^{-1}(W)$
is a $C^{k}$-isomorphism. 
\end{list}

The pair $(U_{\alpha},\boldsymbol{x}_{\alpha})$ with $p\in\boldsymbol{x}_{\alpha}(U_{\alpha})$
is called a \textit{parametrization} (or \textit{system of coordinates})
of $M$ at $p$; and $\boldsymbol{x}_{\alpha}(U_{\alpha})$ is said
to be a \textit{coordinate neighborhood} at $p$.

The set $M$ can be endowed with a topology in a unique way such that
each $\boldsymbol{x}_{\alpha}(U_{\alpha})$ is open, and the $\boldsymbol{x}_{\alpha}$'s
are topological isomorphisms. We note that if $k\geq1$ and two parametrizations
$(U_{\alpha},\boldsymbol{x}_{\alpha})$ and $(U_{\beta},\boldsymbol{x}_{\beta})$
are such that $\boldsymbol{x}_{\alpha}(U_{\alpha})$ and $\boldsymbol{x}_{\beta}(U_{\beta})$
have a non-empty intersection, then from the derivative of $\boldsymbol{x}_{\beta}^{-1}\circ\boldsymbol{x}_{\alpha}$
we have that $X_{\alpha}$ and $X_{\beta}$ are isomorphic.

Two collections $\{(U_{\alpha},\boldsymbol{x}_{\alpha})\}$ and $\{(V_{\beta},\boldsymbol{x}_{\beta})\}$
satisfying (bm1)--(bm3) are said to be \textit{$C^{k}$-compatible}
if their union also satisfies (bm1)--(bm3). It can be verified that
the relation of $C^{k}$-compatibility is an equivalence relation.
An equivalence class of $C^{k}$-compatible collections $\{(U_{\alpha},\boldsymbol{x}_{\alpha})\}$
on $M$ is said to define a \textit{$C^{k}$-differentiable structure}
on $X$.

Now we review the construction of the exponential statistical manifold.
We consider the Musielak--Orlicz space $L^{\Phi_{1}}(p)=L^{\Phi_{1}}(T,\Sigma,p\cdot\mu)$,
where the Orlicz function $\Phi_{1}\colon[0,\infty)\rightarrow[0,\infty)$
is given by $\Phi_{1}(u)=e^{u}-1$, and $p$ is a probability density
in $\mathcal{P}_{\mu}$. The space $L^{\Phi_{1}}(p)$ corresponds
to the set of all functions $u\in L^{0}$ whose \textit{moment-generating
function} $\widehat{u}_{p}(\lambda)=\mathbb{E}_{p}[e^{\lambda u}]$
is finite in a neighborhood of $0$. 

For every function $u\in L^{0}$ we define the \textit{moment-generating
functional} 
\[
M_{p}(u)=\mathbb{E}_{p}[e^{u}],
\]
and the \textit{cumulant-generating functional}
\[
K_{p}(u)=\log M_{p}(u).
\]
Clearly, these functionals are not expected to be finite for every
$u\in L^{0}$. Denote by $\mathcal{K}_{p}$ the interior of the set
of all functions $u\in L^{\Phi_{1}}(p)$ whose moment-generating functional
$M_{p}(u)$ is finite. Equivalently, a function $u\in L^{\Phi_{1}}(p)$
belongs to $\mathcal{K}_{p}$ if and only if $M_{p}(\lambda u)$ is
finite for every $\lambda$ in some neighborhood of $[0,1]$. The
closed subspace of \textit{$p$-centered} random variables
\[
B_{p}=\{u\in L^{\Phi_{1}}(p):\mathbb{E}_{p}[u]=0\}
\]
is taken to be the coordinate Banach space. The \textit{exponential
parametrization} $\boldsymbol{e}_{p}\colon\mathcal{B}_{p}\rightarrow\mathcal{E}_{p}$
maps $\mathcal{B}_{p}=B_{p}\cap\mathcal{K}_{p}$ to the \textit{exponential
family} $\mathcal{E}_{p}=\boldsymbol{e}_{p}(\mathcal{B}_{p})\subseteq\mathcal{P}_{\mu}$,
according to 
\[
\boldsymbol{e}_{p}(u)=e^{u-K_{p}(u)}p,\quad\text{for all }u\in\mathcal{B}_{p}.
\]
$\boldsymbol{e}_{p}$ is a bijection from $\mathcal{B}_{p}$ to its
image $\mathcal{E}_{p}=\boldsymbol{e}_{p}(\mathcal{B}_{p})$, whose
inverse $\boldsymbol{e}_{p}^{-1}\colon\mathcal{E}_{p}\rightarrow\mathcal{B}_{p}$
can be expressed as
\[
\boldsymbol{e}_{p}^{-1}(q)=\log\Bigl(\frac{q}{p}\Bigr)-\mathbb{E}_{p}\Bigl[\log\Bigl(\frac{q}{p}\Bigr)\Bigr],\quad\text{for }q\in\mathcal{E}_{p}.
\]
Since $K_{p}(u)<\infty$ for every $u\in\mathcal{K}_{p}$, we have
that $\boldsymbol{e}_{p}$ can be extended to $\mathcal{K}_{p}$.
The restriction of $\boldsymbol{e}_{p}$ to $\mathcal{B}_{p}$ guarantees
that $\boldsymbol{e}_{p}$ is bijective.

Given two probability densities $p$ and $q$ in the same connected
component of $\mathcal{P}_{\mu}$, the exponential probability families
$\mathcal{E}_{p}$ and $\mathcal{E}_{q}$ coincide, and the exponential
spaces $L^{\Phi_{1}}(p)$ and $L^{\Phi_{1}}(q)$ are isomorphic (see
\cite[Proposition 5]{Pistone:1999}). Hence, $\mathcal{B}_{p}=\boldsymbol{e}_{p}^{-1}(\mathcal{E}_{p}\cap\mathcal{E}_{q})$
and $\mathcal{B}_{q}=\boldsymbol{e}_{q}^{-1}(\mathcal{E}_{p}\cap\mathcal{E}_{q})$.
The transition map $\boldsymbol{e}_{q}^{-1}\circ\boldsymbol{e}_{p}:\mathcal{B}_{p}\rightarrow\mathcal{B}_{q}$,
which can be written as
\[
\boldsymbol{e}_{q}^{-1}\circ\boldsymbol{e}_{p}(u)=u+\log\Bigl(\frac{p}{q}\Bigr)-\mathbb{E}_{q}\Bigl[u+\log\Bigl(\frac{p}{q}\Bigr)\Bigr],\quad\text{for all }u\in\mathcal{B}_{p},
\]
is a $C^{\infty}$-function. Clearly, $\bigcup_{p\in\mathcal{P}_{\mu}}e_{p}(\mathcal{B}_{p})=\mathcal{P}_{\mu}$.
Thus the collection $\{(\mathcal{B}_{p},\boldsymbol{e}_{p})\}_{p\in\mathcal{P}_{\mu}}$
satisfies (bm1)--(bm2). Hence $\mathcal{P}_{\mu}$ is a $C^{\infty}$-Banach
manifold, which is called the \textit{exponential statistical manifold}.

\section{Construction of the $\varphi$-family of probability distributions\label{sec:construction_phi_family}}

The generalization of the exponential family is based on the replacement
of the exponential function by a \textit{$\varphi$-function} $\varphi\colon T\times\overline{\mathbb{R}}\rightarrow[0,\infty]$
that satisfies the following properties, for $\mu$-a.e.\ $t\in T$:
\begin{itemize}
\item [(a1)] $\varphi(t,\cdot)$ is convex and injective,
\item [(a2)] $\varphi(t,-\infty)=0$ and $\varphi(t,\infty)=\infty$,
\item [(a3)] $\varphi(\cdot,u)$ is measurable for all $u\in\mathbb{R}$.
\end{itemize}
In addition, we assume a positive, measurable function $u_{0}\colon T\rightarrow(0,\infty)$
can be found such that, for every measurable function $c\colon T\rightarrow\mathbb{R}$
for which $\varphi(t,c(t))$ is in $\mathcal{P}_{\mu}$, we have that
\begin{itemize}
\item [(a4)] $\varphi(t,c(t)+\lambda u_{0}(t))$ is $\mu$-integrable for
all $\lambda>0$.
\end{itemize}
The choice for $\varphi(t,\cdot)$ injective with image $[0,\infty]$
is justified by the fact that a parametrization of $\mathcal{P}_{\mu}$
maps real-valued functions to positive functions. Moreover, by (a1),
$\varphi(t,\cdot)$ is continuous and strictly increasing. From (a3),
the function $\varphi(t,u(t))$ is measurable if and only if $u\colon T\rightarrow\mathbb{R}$
is measurable. Replacing $\varphi(t,u)$ by $\varphi(t,u_{0}(t)u)$,
a ``new'' function $u_{0}=1$ is obtained satisfying (a4).

\begin{example}
The \textit{Kaniadakis' $\kappa$-exponential} $\exp_{\kappa}\colon\mathbb{R}\rightarrow(0,\infty)$
for $\kappa\in[-1,1]$ is defined as 
\[
\exp_{\kappa}(u)=\begin{cases}
(\kappa u+\sqrt{1+\kappa^{2}u^{2}})^{1/\kappa}, & \text{if }\kappa\neq0,\\
\exp(u), & \text{if }\kappa=0.
\end{cases}
\]
The inverse of $\exp_{\kappa}$ is the \textit{Kaniadakis' $\kappa$-logarithm}
\[
\ln_{\kappa}(u)=\begin{cases}
\dfrac{u^{\kappa}-u^{-\kappa}}{2\kappa}, & \text{if }\kappa\neq0,\\
\ln(u), & \text{if }\kappa=0.
\end{cases}
\]
Some algebraic properties of the ordinary exponential and logarithm
functions are preserved:
\[
\exp_{\kappa}(u)\exp_{\kappa}(-u)=1,\qquad\ln_{\kappa}(u)+\ln_{\kappa}(u^{-1})=0.
\]
For a measurable function $\kappa\colon T\rightarrow[-1,1],$ we define
the \textit{variable $\kappa$-exponential} $\exp_{\kappa}\colon T\times\mathbb{R}\rightarrow(0,\infty)$
as 
\[
\exp_{\kappa}(t,u)=\exp_{\kappa(t)}(u),
\]
whose inverse is called the \textit{variable $\kappa$-logarithm}:
\[
\ln_{\kappa}(t,u)=\ln_{\kappa(t)}(u).
\]
Assuming that $\kappa_{-}=\operatorname{ess\, inf}|\kappa(t)|>0$,
the variable $\kappa$-exponential $\exp_{\kappa}$ satisfies (a1)--(a4).
The verification of (a1)--(a3) is easy. Moreover, we notice that $\exp_{\kappa}(t,\cdot)$
is strictly convex. We can write for $\alpha\geq1$
\begin{align*}
\exp_{\kappa}(t,\alpha u) & =(\kappa(t)\alpha u+\alpha\sqrt{1/\alpha^{2}+\kappa(t)^{2}u^{2}})^{1/\kappa}\\
 & \leq\alpha^{1/|\kappa|}(\kappa(t)u+\sqrt{1+\kappa(t)^{2}u^{2}})^{1/\kappa}\\
 & \leq\alpha^{1/\kappa_{-}}\exp_{\kappa}(t,u).
\end{align*}
By the convexity of $\exp_{\kappa}(t,\cdot)$, we obtain for any $\lambda\in(0,1)$
\begin{align*}
\exp_{\kappa}(t,c+u) & \leq\lambda\exp_{\kappa}(t,\lambda^{-1}c)+(1-\lambda)\exp_{\kappa}(t,(1-\lambda)^{-1}u)\\
 & \leq\lambda^{1-1/\kappa_{-}}\exp_{\kappa}(t,c)+(1-\lambda)^{1-1/\kappa_{-}}\exp_{\kappa}(t,u).
\end{align*}
Thus any positive function $u_{0}$ such that $\mathbb{E}[\exp_{\kappa}(u_{0})]<\infty$
satisfies (a4).
\end{example}

Let $c\colon T\rightarrow\mathbb{R}$ be a measurable function such
that $\varphi(t,c(t))$ is $\mu$-integrable. We define the Musielak--Orlicz
function 
\[
\Phi(t,u)=\varphi(t,c(t)+u)-\varphi(t,c(t)).
\]
and denote $L^{\Phi}$, $\tilde{L}^{\Phi}$ and $E^{\Phi}$ by $L_{c}^{\varphi}$,
$\tilde{L}_{c}^{\varphi}$ and $E_{c}^{\varphi}$, respectively. Since
$\varphi(t,c(t))$ is $\mu$-integrable, the Musielak--Orlicz space
$L_{c}^{\varphi}$ corresponds to the set of all functions $u\in L^{0}$
for which $\varphi(t,c(t)+\lambda u(t))$ is $\mu$-integrable for
every $\lambda$ contained in some neighborhood of $0$. 

Let $\mathcal{K}_{c}^{\varphi}$ be the set of all functions $u\in L_{c}^{\varphi}$
such that $\varphi(t,c(t)+\lambda u(t))$ is $\mu$-integrable for
every $\lambda$ in a neighborhood of $[0,1]$. Denote by $\boldsymbol{\varphi}$
the operator acting on the set of real-valued functions $u\colon T\rightarrow\mathbb{R}$
given by $\boldsymbol{\varphi}(u)(t)=\varphi(t,u(t))$. For each probability
density $p\in\mathcal{P}_{\mu}$, we can take a measurable function
$c\colon T\rightarrow\mathbb{R}$ such that $p=\boldsymbol{\varphi}(c)$.
The first import result in the construction of the $\varphi$-family
is given below.

\begin{lem}
The set $\mathcal{K}_{c}^{\varphi}$ is open in $L_{c}^{\varphi}$.\end{lem}
\begin{proof}
Take any $u\in\mathcal{K}_{c}^{\varphi}$. We can find $\varepsilon\in(0,1)$
such that $\mathbb{E}[\boldsymbol{\varphi}(c+\alpha u)]<\infty$ for
every $\alpha\in[-\varepsilon,1+\varepsilon]$. Let $\delta=[\frac{2}{\varepsilon}(1+\varepsilon)(1+\frac{\varepsilon}{2})]^{-1}$.
For any function $v\in L_{c}^{\varphi}$ in the open ball $B_{\delta}=\{w\in L_{c}^{\varphi}:\Vert w\Vert_{\Phi}<\delta\}$,
we have $I_{\Phi}(\frac{v}{\delta})\leq1$. Thus $\mathbb{E}[\boldsymbol{\varphi}(c+\frac{1}{\delta}|v|)]\leq2$.
Taking any $\alpha\in(0,1+\frac{\varepsilon}{2})$, we denote $\lambda=\frac{\alpha}{1+\varepsilon}$.
In virtue of
\[
\frac{\alpha}{1-\lambda}=\frac{\alpha}{1-\frac{\alpha}{1+\varepsilon}}\leq\frac{1+\frac{\varepsilon}{2}}{1-\frac{1+\frac{\varepsilon}{2}}{1+\varepsilon}}=\frac{2}{\varepsilon}(1+\varepsilon)(1+\frac{\varepsilon}{2})=\frac{1}{\delta},
\]
it follows that
\begin{align}
\boldsymbol{\varphi}(c+\alpha(u+v)) & =\boldsymbol{\varphi}(\lambda(c+\tfrac{\alpha}{\lambda}u)+(1-\lambda)(c+\tfrac{\alpha}{1-\lambda}v))\nonumber \\
 & \leq\lambda\boldsymbol{\varphi}(c+\tfrac{\alpha}{\lambda}u)+(1-\lambda)\boldsymbol{\varphi}(c+\tfrac{\alpha}{1-\lambda}v)\nonumber \\
 & \leq\lambda\boldsymbol{\varphi}(c+(1+\varepsilon)u)+(1-\lambda)\boldsymbol{\varphi}(c+\tfrac{1}{\delta}|v|).\label{eq:open_1}
\end{align}
For $\alpha\in(-\frac{\varepsilon}{2},0)$, we can write
\begin{align}
\boldsymbol{\varphi}(c+\alpha(u+v)) & \leq\tfrac{1}{2}\boldsymbol{\varphi}(c+2\alpha u)+\tfrac{1}{2}\boldsymbol{\varphi}(c+2\alpha v)\nonumber \\
 & \leq\tfrac{1}{2}\boldsymbol{\varphi}(c+2\alpha u)+\tfrac{1}{2}\boldsymbol{\varphi}(c+|v|).\label{eq:open_2}
\end{align}
By \prettyref{eq:open_1} and \prettyref{eq:open_2}, we get $\mathbb{E}[\boldsymbol{\varphi}(c+\alpha(u+v))]<\infty$,
for any $\alpha\in(-\frac{\varepsilon}{2},1+\frac{\varepsilon}{2})$.
Hence the ball of radius $\delta$ centered at $u$ is contained in
$\mathcal{K}_{c}^{\varphi}$. Therefore, the set $\mathcal{K}_{c}^{\varphi}$
is open.
\end{proof}

Clearly, for $u\in\mathcal{K}_{c}^{\varphi}$ the function $\boldsymbol{\varphi}(c+u)$
is not necessarily in $\mathcal{P}_{\mu}$. The \textit{normalizing
function} $\psi\colon\mathcal{K}_{c}^{\varphi}\rightarrow\mathbb{R}$
is introduced in order to make the density 
\[
\boldsymbol{\varphi}(c+u-\psi(u)u_{0})
\]
contained in $\mathcal{P}_{\mu}$, for any $u\in\mathcal{K}_{c}^{\varphi}$.
We have to find the functions for which the normalizing function exists.
For a function $u\in L_{c}^{\varphi}$, suppose that $\boldsymbol{\varphi}(c+u-\alpha u_{0})$
is $\mu$-integrable for some $\alpha\in\mathbb{R}$. Then $u$ is
in the closure of the set $\mathcal{K}_{c}^{\varphi}$. Indeed, for
any $\lambda\in(0,1)$,
\begin{align*}
\boldsymbol{\varphi}(c+\lambda u) & =\boldsymbol{\varphi}(\lambda(c+u-\alpha u_{0})+(1-\lambda)(c+\tfrac{\lambda}{1-\lambda}\alpha u_{0}))\\
 & \leq\lambda\boldsymbol{\varphi}(c+u-\alpha u_{0})+(1-\lambda)\boldsymbol{\varphi}(c+\tfrac{\lambda}{1-\lambda}\alpha u_{0}).
\end{align*}
Since the function $u_{0}$ satisfies (a4), we obtain that $\boldsymbol{\varphi}(c+\lambda u)$
is $\mu$-integrable. Hence the maximal, open domain of $\psi$ is
contained in $\mathcal{K}_{c}^{\varphi}$.

\begin{prop}
\label{pro:ExUn} If the function $u$ is in $\mathcal{K}_{c}^{\varphi}$,
then there exists a unique $\psi(u)\in\mathbb{R}$ for which $\boldsymbol{\varphi}(c+u-\psi(u)u_{0})$
is a probability density in $\mathcal{P}_{\mu}$.\end{prop}
\begin{proof}
We will show that if the function $u$ is in $\mathcal{K}_{c}^{\varphi}$,
then $\boldsymbol{\varphi}(c+u+\alpha u_{0})$ is $\mu$-integrable
for every $\alpha\in\mathbb{R}$. Since $u$ is in $\mathcal{K}_{c}^{\varphi}$,
we can find $\varepsilon>0$ such that $\boldsymbol{\varphi}(c+(1+\varepsilon)u)$
is $\mu$-integrable. Taking $\lambda=\frac{1}{1+\varepsilon}$, we
can write
\begin{align*}
\boldsymbol{\varphi}(c+u+\alpha u_{0}) & =\boldsymbol{\varphi}(\lambda(c+\tfrac{1}{\lambda}u)+(1-\lambda)(c+\tfrac{1}{1-\lambda}\alpha u_{0}))\\
 & \leq\lambda\boldsymbol{\varphi}(c+\tfrac{1}{\lambda}u)+(1-\lambda)\boldsymbol{\varphi}(c+\tfrac{1}{1-\lambda}\alpha u_{0}).
\end{align*}
Thus $\boldsymbol{\varphi}(c+u+\alpha u_{0})$ is $\mu$-integrable.
By the Dominated Convergence Theorem, the map $\alpha\mapsto J(\alpha)=\mathbb{E}[\boldsymbol{\varphi}(c+u+\alpha u_{0})]$
is continuous, tends to $0$ as $\alpha\rightarrow-\infty$, and goes
to infinity as $\alpha\rightarrow\infty$. Since $\varphi(t,\cdot)$
is strictly increasing, it follows that $J(\alpha)$ is also strictly
increasing. Therefore, there exists a unique $\psi(u)\in\mathbb{R}$
for which $\boldsymbol{\varphi}(c+u-\psi(u)u_{0})$ is a probability
density in $\mathcal{P}_{\mu}$.
\end{proof}

The function $\psi\colon\mathcal{K}_{c}^{\varphi}\rightarrow\mathbb{R}$
can take both positive and negative values. However, if the domain
of $\psi$ is restricted to a subspace of $L_{c}^{\varphi}$, its
image will be contained in $[0,\infty)$. We denote by $\boldsymbol{\varphi}_{+}'$
the operator acting on the set of real-valued functions $u\colon T\rightarrow\mathbb{R}$
given by $\boldsymbol{\varphi}_{+}'(u)(t)=\varphi_{+}'(t,u(t))$,
where $\varphi_{+}'(t,\cdot)$ is the right-derivative of $\varphi(t,\cdot)$.
Define the closed subspace
\[
B_{c}^{\varphi}=\{u\in L_{c}^{\varphi}:\mathbb{E}[u\boldsymbol{\varphi}_{+}'(c)]=0\}.
\]
and let $\mathcal{B}_{c}^{\varphi}=B_{c}^{\varphi}\cap\mathcal{K}_{c}^{\varphi}$.
By the convexity of $\varphi(t,\cdot)$, we have 
\[
u\varphi_{+}'(t,c(t))\leq\varphi(t,c(t)+u)-\varphi(t,c(t)),\text{ for all }u\in\mathbb{R}.
\]
Hence, for any $u\in\mathcal{B}_{c}^{\varphi}$, we get
\[
1=\mathbb{E}[u\boldsymbol{\varphi}_{+}'(c)]+\mathbb{E}[\boldsymbol{\varphi}(c)]\leq\mathbb{E}[\boldsymbol{\varphi}(c+u)]<\infty.
\]
Thus it follows that $\psi(u)\geq0$ in order to obtain that $\boldsymbol{\varphi}(c+u-\psi(u)u_{0})$
is in $\mathcal{P}_{\mu}$.

For each measurable function $c\colon T\rightarrow\mathbb{R}$ such
that $p=\boldsymbol{\varphi}(c)$ is the probability density in $\mathcal{P}_{\mu}$,
we associate a parametrization $\boldsymbol{\varphi}_{c}\colon\mathcal{B}_{c}^{\varphi}\rightarrow\mathcal{F}_{c}^{\varphi}$
that maps any function $u$ in $\mathcal{B}_{c}^{\varphi}$ to a probability
density in $\mathcal{F}_{c}^{\varphi}=\varphi_{c}(\mathcal{B}_{c}^{\varphi})\subseteq\mathcal{P}_{\mu}$
according to 
\[
\boldsymbol{\varphi}_{c}(u)=\boldsymbol{\varphi}(c+u-\psi(u)u_{0}).
\]
Clearly, we have $\mathcal{P}_{\mu}=\bigcup\{\mathcal{F}_{c}^{\varphi}:\boldsymbol{\varphi}(c)\in\mathcal{P}_{\mu}\}$.
Moreover, the map $\boldsymbol{\varphi}_{c}$ is a bijection from
$\mathcal{B}_{c}^{\varphi}$ to $\mathcal{F}_{c}^{\varphi}$. If the
functions $u,v\in\mathcal{B}_{c}^{\varphi}$ are such that $\boldsymbol{\varphi}_{c}(u)=\boldsymbol{\varphi}_{c}(v)$,
then the difference $u-v=(\psi(u)-\psi(v))u_{0}$ is in $B_{c}^{\varphi}$.
Consequently, $\psi(u)=\psi(v)$ and then $u=v$.

Suppose that the measurable functions $c_{1},c_{2}\colon T\rightarrow\mathbb{R}$
are such that $p_{1}=\boldsymbol{\varphi}(c_{1})$ and $p_{2}=\boldsymbol{\varphi}(c_{2})$
belong to $\mathcal{P}_{\mu}$. The parametrizations $\boldsymbol{\varphi}_{c_{1}}\colon\mathcal{B}_{c_{1}}^{\varphi}\rightarrow\mathcal{F}_{c_{1}}^{\varphi}$
and $\boldsymbol{\varphi}_{c_{2}}\colon\mathcal{B}_{c_{2}}^{\varphi}\rightarrow\mathcal{F}_{c_{2}}^{\varphi}$
related to these functions have transition map 
\[
\boldsymbol{\varphi}_{c_{2}}^{-1}\circ\boldsymbol{\varphi}_{c_{1}}\colon\boldsymbol{\varphi}_{c_{1}}^{-1}(\mathcal{F}_{c_{1}}^{\varphi}\cap\mathcal{F}_{c_{2}}^{\varphi})\rightarrow\boldsymbol{\varphi}_{c_{2}}^{-1}(\mathcal{F}_{c_{1}}^{\varphi}\cap\mathcal{F}_{c_{2}}^{\varphi}).
\]
Let $\psi_{1}\colon\mathcal{B}_{c_{1}}^{\varphi}\rightarrow[0,\infty)$
and $\psi_{2}\colon\mathcal{B}_{c_{2}}^{\varphi}\rightarrow[0,\infty)$
be the normalizing functions associated to $c_{1}$ and $c_{2}$,
respectively. Assume that the functions $u\in\mathcal{B}_{c_{1}}^{\varphi}$
and $v\in\mathcal{B}_{c_{2}}^{\varphi}$ are such that $\boldsymbol{\varphi}_{c_{1}}(u)=\boldsymbol{\varphi}_{c_{2}}(v)\in\mathcal{F}_{c_{1}}^{\varphi}\cap\mathcal{F}_{c_{2}}^{\varphi}$.
Then we can write
\[
v=c_{1}-c_{2}+u-(\psi_{1}(u)-\psi_{2}(v))u_{0}.
\]
Since the function $v$ is in $B_{c_{2}}^{\varphi}$, if we multiply
this equation by $\boldsymbol{\varphi}_{+}'(c_{2})$ and integrate
with respect to the measure $\mu$, we obtain
\[
0=\mathbb{E}[(c_{1}-c_{2}+u)\boldsymbol{\varphi}_{+}'(c_{2})]-(\psi_{1}(u)-\psi_{2}(v))\mathbb{E}[u_{0}\boldsymbol{\varphi}_{+}'(c_{2})].
\]
Thus the transition map $\boldsymbol{\varphi}_{c_{2}}^{-1}\circ\boldsymbol{\varphi}_{c_{1}}$
can be expressed as 
\begin{equation}
\boldsymbol{\varphi}_{c_{2}}^{-1}\circ\boldsymbol{\varphi}_{c_{1}}(w)=c_{1}-c_{2}+w-\frac{\mathbb{E}[(c_{1}-c_{2}+w)\boldsymbol{\varphi}_{+}'(c_{2})]}{\mathbb{E}[u_{0}\boldsymbol{\varphi}_{+}'(c_{2})]}u_{0},\label{eq:trans_map}
\end{equation}
for every $w\in\boldsymbol{\varphi}_{c_{1}}^{-1}(\mathcal{F}_{c_{1}}^{\varphi}\cap\mathcal{F}_{c_{2}}^{\varphi})$.
Clearly, this transition map will be of class $C^{\infty}$ if we
show that the functions $w$ and $c_{1}-c_{2}$ are in $L_{c_{2}}^{\varphi}$,
and the spaces $L_{c_{1}}^{\varphi}$ and $L_{c_{2}}^{\varphi}$ have
equivalent norms. It is not hard to verify that if two Musielak--Orlicz
spaces are equal as sets, then their norms are equivalent (see \cite[Theorem 8.5]{Musielak:1983}).
We make use of the following:

\begin{prop}
\label{pro:incl_ct_c} Assume that the measurable functions $\widetilde{c},c\colon T\rightarrow\mathbb{R}$
satisfy $\mathbb{E}[\varphi(t,\widetilde{c}(t))]<\infty$ and $\mathbb{E}[\varphi(t,c(t))]<\infty$.
Then $L_{\widetilde{c}}^{\varphi}\subseteq L_{c}^{\varphi}$ if and
only if $\widetilde{c}-c\in L_{c}^{\varphi}$.\end{prop}
\begin{proof}
Suppose that $\widetilde{c}-c$ is not in $L_{c}^{\varphi}$. Let
$A=\{t\in T:\widetilde{c}(t)<c(t)\}$. For $\lambda\in[0,1]$, we
have 
\begin{align*}
\mathbb{E}[\boldsymbol{\varphi}(c+\lambda(\widetilde{c}-c))] & =\mathbb{E}[\boldsymbol{\varphi}(c+\lambda(\widetilde{c}-c))\boldsymbol{1}_{T\setminus A}]+\mathbb{E}[\boldsymbol{\varphi}(c+\lambda(\widetilde{c}-c))\boldsymbol{1}_{A}]\\
 & \leq\mathbb{E}[\boldsymbol{\varphi}(c+(\widetilde{c}-c))\boldsymbol{1}_{T\setminus A}]+\mathbb{E}[\boldsymbol{\varphi}(c)\boldsymbol{1}_{A}]\\
 & \leq\mathbb{E}[\boldsymbol{\varphi}(\widetilde{c})]+\mathbb{E}[\boldsymbol{\varphi}(c)]<\infty.
\end{align*}
Since $\widetilde{c}-c\notin L_{c}^{\varphi}$, for any $\lambda>0$,
there holds $\mathbb{E}[\boldsymbol{\varphi}(c-\lambda(\widetilde{c}-c))]=\infty$.
From 
\begin{align*}
\mathbb{E}[\boldsymbol{\varphi}(c-\lambda(\widetilde{c}-c))] & =\mathbb{E}[\boldsymbol{\varphi}(c-\lambda(\widetilde{c}-c))\boldsymbol{1}_{T\setminus A}]+\mathbb{E}[\boldsymbol{\varphi}(c-\lambda(\widetilde{c}-c))\boldsymbol{1}_{A}]\\
 & \leq\mathbb{E}[\boldsymbol{\varphi}(c+\lambda(c-\widetilde{c}))\boldsymbol{1}_{A}],
\end{align*}
we obtain that $(c-\widetilde{c})\boldsymbol{1}_{A}$ does not belong
to $L_{c}^{\varphi}$. Clearly, $(c-\widetilde{c})\boldsymbol{1}_{A}\in L_{\widetilde{c}}^{\varphi}$.
Consequently, $L_{\widetilde{c}}^{\varphi}$ is not contained in $L_{c}^{\varphi}$.

Conversely, assume $\widetilde{c}-c\in L_{c}^{\varphi}$. Let $w$
be any function in $L_{\widetilde{c}}^{\varphi}$. We can find $\varepsilon>0$
such that $\mathbb{E}[\boldsymbol{\varphi}(\widetilde{c}+\lambda w)]<\infty$,
for every $\lambda\in(-\varepsilon,\varepsilon)$. Consider the convex
function
\[
g(\alpha,\lambda)=\mathbb{E}[\boldsymbol{\varphi}(c+\alpha(\widetilde{c}-c)+\lambda w)].
\]
This function is finite for $\lambda=0$ and $\alpha$ in the interval
$(-\eta,1]$, for some $\eta>0$. Moreover, $g(1,\lambda)$ is finite
for every $\lambda\in(-\varepsilon,\varepsilon)$. By the convexity
of $g$, we have that $g$ is finite in the convex hull of the set
$1\times(-\varepsilon,\varepsilon)\cup(-\eta,1]\times0$. We obtain
that $g(0,\lambda)$ is finite for every $\lambda$ in some neighborhood
of $0$. Consequently, $w\in L_{c}^{\varphi}$. Since $w\in L_{c}^{\varphi}$
is arbitrary, the inclusion $L_{\widetilde{c}}^{\varphi}\subseteq L_{c}^{\varphi}$
follows.
\end{proof}

\begin{lem}
\label{lem:eqLP} If the function $u$ is in $\mathcal{K}_{c}^{\varphi}$
and we denote $\widetilde{c}=c+u-\psi(u)u_{0}$, then the spaces $L_{c}^{\varphi}$
and $L_{\widetilde{c}}^{\varphi}$ are equal as sets.\end{lem}
\begin{proof}
The inclusion $L_{\widetilde{c}}^{\varphi}\subseteq L_{c}^{\varphi}$
follows from \prettyref{pro:incl_ct_c}. Since $u\in\mathcal{K}_{c}^{\varphi}$,
we have 
\[
\mathbb{E}[\boldsymbol{\varphi}(\widetilde{c}+\lambda u)]\leq\mathbb{E}[\boldsymbol{\varphi}(c+(1+\lambda)u)]<\infty,
\]
for every $\lambda$ in a neighborhood of $0$. Thus $c-\widetilde{c}=-u+\psi(u)u_{0}$
belongs to $L_{\widetilde{c}}^{\varphi}$. From \prettyref{pro:incl_ct_c},
we obtain $L_{\widetilde{c}}^{\varphi}\subseteq L_{c}^{\varphi}$.
\end{proof}

By \prettyref{lem:eqLP}, if we denote $c_{1}+u-\psi_{1}(u)u_{0}=\widetilde{c}=c_{2}+v-\psi_{2}(v)u_{0}$,
we have that the spaces $L_{c_{1}}^{\varphi}$, $L_{\widetilde{c}}^{\varphi}$
and $L_{c_{2}}^{\varphi}$ are equal as sets. In \prettyref{eq:trans_map},
the function $w$ is in $L_{c_{2}}^{\varphi}$ and consequently $c_{1}-c_{2}$
is in $L_{c_{2}}^{\varphi}$. Therefore, the transition map $\boldsymbol{\varphi}_{c_{2}}^{-1}\circ\boldsymbol{\varphi}_{c_{1}}$
is of class $C^{\infty}$.

Since $\boldsymbol{\varphi}_{c_{2}}^{-1}\circ\boldsymbol{\varphi}_{c_{1}}$
is of class $C^{\infty}$, the set $\boldsymbol{\varphi}_{c_{1}}^{-1}(\mathcal{F}_{c_{1}}^{\varphi}\cap\mathcal{F}_{c_{2}}^{\varphi})$
is open $B_{c_{1}}^{\varphi}$. The $\varphi$-families $\mathcal{F}_{c}^{\varphi}$
are maximal in the sense that if two $\varphi$-families $\mathcal{F}_{c_{1}}^{\varphi}$
and $\mathcal{F}_{c_{2}}^{\varphi}$ have non-empty intersection,
then they coincide.

\begin{lem}
\label{lem:eqModel} For a function $u$ in $\mathcal{B}_{c}^{\varphi}$,
denote $\widetilde{c}=c+u-\psi(u)u_{0}$. Then $\mathcal{F}_{c}^{\varphi}=\mathcal{F}_{\widetilde{c}}^{\varphi}$.\end{lem}
\begin{proof}
Let $v$ be a function in $\mathcal{B}_{c}^{\varphi}$. Then there
exists $\varepsilon>0$ such that, for every $\lambda\in(-\varepsilon,1+\varepsilon)$,
the function $\boldsymbol{\varphi}(c+\lambda v+(1-\lambda)u)$ is
$\mu$-integrable. Consequently, $\varphi(\widetilde{c}+\lambda(v-u))$
is $\mu$-integrable for all $\lambda\in(-\varepsilon,1+\varepsilon)$.
Thus the difference $v-u$ is in $\mathcal{K}_{\widetilde{c}}^{\varphi}$
and 
\begin{equation}
w=v-u-\frac{\mathbb{E}[(v-u)\boldsymbol{\varphi}_{+}'(\widetilde{c})]}{\mathbb{E}[u_{0}\boldsymbol{\varphi}_{+}'(\widetilde{c})]}u_{0}\label{eq:proj_dif}
\end{equation}
belongs to $\mathcal{B}_{\widetilde{c}}^{\varphi}$. Let $\widetilde{\psi}\colon\mathcal{B}_{\widetilde{c}}^{\varphi}\rightarrow[0,\infty)$
be the normalizing function associated to $\widetilde{c}$. Then the
probability density $\boldsymbol{\varphi}(\widetilde{c}+w-\widetilde{\psi}(w)u_{0})$
is in $\mathcal{F}_{\widetilde{c}}^{\varphi}$. This probability density
can be expressed as $\boldsymbol{\varphi}(c+v-ku_{0})$ for a constant
$k$. According to \prettyref{pro:ExUn}, there exists a unique $\psi(u)\in\mathbb{R}$
such that the probability density $\boldsymbol{\varphi}(c+v-\psi(v)u_{0})$
is in $\mathcal{F}_{c}^{\varphi}$. Therefore, $\mathcal{F}_{c}^{\varphi}\subseteq\mathcal{F}_{\widetilde{c}}^{\varphi}$.

Using the same arguments as in the previous paragraph, we obtain that
$c=\widetilde{c}+w-\widetilde{\psi}(w)u_{0}$, where the function
$w\in\mathcal{B}_{\widetilde{c}}^{\varphi}$ is given in \prettyref{eq:proj_dif}
with $v=0$. Thus $\mathcal{F}_{\widetilde{c}}^{\varphi}\subseteq\mathcal{F}_{c}^{\varphi}$.
\end{proof}

By \prettyref{lem:eqModel}, if we denote $c_{1}+u-\psi_{1}(u)u_{0}=\widetilde{c}=c_{2}+v-\psi_{2}(v)u_{0}$,
then we have the equality $\mathcal{F}_{c_{1}}^{\varphi}=\mathcal{F}_{\widetilde{c}}^{\varphi}=\mathcal{F}_{c_{2}}^{\varphi}$.

The results obtained in these lemmas are summarized in the next Proposition.

\begin{prop}
Let $c_{1},c_{2}\colon T\rightarrow\mathbb{R}$ be measurable functions
such that the probability densities $p_{1}=\boldsymbol{\varphi}(c_{1})$
and $p_{2}=\boldsymbol{\varphi}(c_{2})$ are in $\mathcal{P}_{\mu}$.
Suppose $\mathcal{F}_{c_{1}}^{\varphi}\cap\mathcal{F}_{c_{2}}^{\varphi}\neq\emptyset$.
Then the Musielak--Orlicz spaces $L_{c_{1}}^{\varphi}$ and $L_{c_{2}}^{\varphi}$
are equal as sets, and have equivalent norms. Moreover, $\mathcal{F}_{c_{1}}^{\varphi}=\mathcal{F}_{c_{2}}^{\varphi}$\textup{.}
\end{prop}

Thus we can state:

\begin{prop}
The collection $\{(\mathcal{B}_{c}^{\varphi},\boldsymbol{\varphi}_{c})\}_{\boldsymbol{\varphi}(c)\in\mathcal{P}_{\mu}}$
satisfies (bm1)--(bm2), equipping $\mathcal{P}_{\mu}$ with a $C^{\infty}$-differentiable
structure.
\end{prop}

\section{Divergence\label{sec:divergence}}

In this section we define the divergence between two probability distributions.
The entities found in Information Geometry \cite{Amari:2000,Murray:1993},
like the Fisher information, connections, geodesics, etc., are all
derived from the divergence taken in the considered family. The divergence
we will found is the Bregman divergence \cite{Bregman:1967} associated
to the normalizing function $\psi\colon\mathcal{K}_{c}^{\varphi}\rightarrow[0,\infty)$.
We show that our divergence does not depend on the parametrization
of the $\varphi$-family $\mathcal{F}_{c}^{\varphi}$.

Let $S$ be a convex subset of a Banach space $X$. Given a convex
function $f\colon S\rightarrow\mathbb{R}$, the \textit{Bregman divergence}
$B_{f}\colon S\times S\rightarrow[0,\infty)$ is defined as 
\[
B_{f}(y,x)=f(y)-f(x)-\partial_{+}f(x)(y-x),
\]
for all $x,y\in S$, where $\partial_{+}f(x)(h)=\lim_{t\downarrow0}(f(x+th)-f(x))/t$
denotes the \textit{right-directional derivative} of $f$ at $x$
in the direction of $h$. The right-directional derivative $\partial_{+}f(x)(h)$
exists and defines a sublinear functional. If the function $f$ is
strictly convex, the divergence satisfies $B_{f}(y,x)=0$ if and only
if $x=y$. 

Let $X$ and $Y$ be Banach spaces, and $U\subseteq X$ be an open
set. A function $f\colon U\rightarrow Y$ is said to be \textit{Gâteaux-differentiable}
at $x_{0}\in U$ if there exists a bounded linear map $A\colon X\rightarrow Y$
such that 
\[
\lim_{t\rightarrow0}\frac{1}{t}\Vert f(x_{0}+th)-f(x_{0})-Ah\Vert=0,
\]
for every $h\in X$. The \textit{Gâteaux derivative} of $f$ at $x_{0}$
is denoted by $A=\partial f(x_{0})$. If the limit above can be taken
uniformly for every $h\in X$ such that $\Vert h\Vert\leq1$, then
the function $f$ is said to be \textit{Fréchet-differentiable} at
$x_{0}$. The \textit{Fréchet derivative} of $f$ at $x_{0}$ is denoted
by $A=Df(x_{0})$.

Now we verify that $\psi\colon\mathcal{K}_{c}^{\varphi}\rightarrow\mathbb{R}$
is a convex function. Take any $u,v\in\mathcal{K}_{c}^{\varphi}$
such that $u\neq v$. Clearly, the function $\lambda u+(1-\lambda)v$
is in $\mathcal{K}_{c}^{\varphi}$, for any $\lambda\in(0,1)$. By
the convexity of $\varphi(t,\cdot)$, we can write
\begin{multline*}
\mathbb{E}[\boldsymbol{\varphi}(c+\lambda u+(1-\lambda)v-\lambda\psi(u)u_{0}-(1-\lambda)\psi(v)u_{0})]\\
\leq\lambda\mathbb{E}[\boldsymbol{\varphi}(c+u-\psi(u)u_{0})]+(1-\lambda)\mathbb{E}[\boldsymbol{\varphi}(c+v-\psi(v)u_{0})]=1.
\end{multline*}
Since $\boldsymbol{\varphi}(c+\lambda u+(1-\lambda)v-\psi(\lambda u+(1-\lambda)v)u_{0})$
has $\mu$-integral equal to $1$, we can conclude that the following
inequality holds:
\[
\psi(\lambda u+(1-\lambda)v)\leq\lambda\psi(u)+(1-\lambda)\psi(v).
\]
So we can define the Bregman divergence $B_{\psi}$ from to the normalizing
function $\psi$.

The Bregman divergence $B_{\psi}\colon\mathcal{B}_{c}^{\varphi}\times\mathcal{B}_{c}^{\varphi}\rightarrow[0,\infty)$
associated to the normalizing function $\psi\colon\mathcal{B}_{c}^{\varphi}\rightarrow[0,\infty)$
is given by
\[
B_{\psi}(v,u)=\psi(v)-\psi(u)-\partial_{+}\psi(u)(v-u).
\]
Then we define the divergence $D_{\psi}\colon\mathcal{B}_{c}^{\varphi}\times\mathcal{B}_{c}^{\varphi}\rightarrow[0,\infty)$
related to the $\varphi$-family $\mathcal{F}_{c}^{\varphi}$ as
\[
D_{\psi}(u,v)=B_{\psi}(v,u).
\]
The entries of $B_{\psi}$ are inverted in order that $D_{\psi}$
corresponds in some way to the \textit{Kullback--Leibler divergence}
$D_{\mathrm{KL}}(p,q)=\mathbb{E}[p\log(\frac{p}{q})]$. Assuming that
$\varphi(t,\cdot)$ is continuously differentiable, we will find an
expression for $\partial\psi(u)$.

\begin{lem}
\label{lem:func_phi_prime} Assume that $\varphi(t,\cdot)$ is continuously
differentiable. For any $u\in\mathcal{K}_{c}^{\varphi}$, the linear
functional $f_{u}\colon L_{c}^{\varphi}\rightarrow\mathbb{R}$ given
by $f_{u}(v)=\mathbb{E}[v\boldsymbol{\varphi}'(c+u)]$ is bounded.\end{lem}
\begin{proof}
Every function $v\in L_{c}^{\varphi}$ with norm $\Vert v\Vert_{\Phi,0}\leq1$
satisfies $I_{\Phi}(v)\leq\Vert v\Vert_{\Phi,0}$. Then we obtain
\[
\mathbb{E}[\boldsymbol{\varphi}(c+|v|)]=I_{\Phi}(v)+\mathbb{E}[\boldsymbol{\varphi}(c)]\leq2.
\]
Since $u\in\mathcal{K}_{c}^{\varphi}$, we can find $\lambda\in(0,1)$
such that $\mathbb{E}[\boldsymbol{\varphi}(c+\tfrac{1}{\lambda}u)]<\infty$.
We can write
\begin{align*}
(1-\lambda)\mathbb{E}[|v|\boldsymbol{\varphi}'(c+u)] & \leq\mathbb{E}[\boldsymbol{\varphi}(c+u+(1-\lambda)|v|)]-\mathbb{E}[\boldsymbol{\varphi}(c+u)]\\
 & =\mathbb{E}[\boldsymbol{\varphi}(\lambda(c+\tfrac{1}{\lambda}u)+(1-\lambda)(c+|v|))]-\mathbb{E}[\boldsymbol{\varphi}(c+u)]\\
 & \leq\lambda\mathbb{E}[\boldsymbol{\varphi}(c+\tfrac{1}{\lambda}u)]+(1-\lambda)\mathbb{E}[\boldsymbol{\varphi}(c+|v|)]-\mathbb{E}[\boldsymbol{\varphi}(c+u)].
\end{align*}
Thus the absolute value of $f_{u}(v)=\mathbb{E}[v\boldsymbol{\varphi}'(c+u)]$
is bounded by some constant for $\Vert v\Vert_{\Phi,0}\leq1$.
\end{proof}

\begin{lem}
Assume that $\varphi(t,\cdot)$ is continuously differentiable. Then
the normalizing function $\psi\colon\mathcal{K}_{c}^{\varphi}\rightarrow\mathbb{R}$
is Gâteaux-differentiable and \textup{
\begin{equation}
\partial\psi(u)v=\frac{\mathbb{E}[v\boldsymbol{\varphi}'(c+u-\psi(u)u_{0})]}{\mathbb{E}[u_{0}\boldsymbol{\varphi}'(c+u-\psi(u)u_{0})]}.\label{eq:partial_psi}
\end{equation}
}\end{lem}
\begin{proof}
According to \prettyref{lem:func_phi_prime}, the expression in \prettyref{eq:partial_psi}
defines a bounded linear functional. Fix functions $u\in\mathcal{K}_{c}^{\varphi}$
and $v\in L_{c}^{\varphi}$. In virtue of \prettyref{pro:incl_ct_c},
we can find $\varepsilon>0$ such that $\mathbb{E}[\boldsymbol{\varphi}(c+u+\lambda|v|)]<\infty$,
for every $\lambda\in[-\varepsilon,\varepsilon]$. Define
\[
g(\lambda,k)=\mathbb{E}[\boldsymbol{\varphi}(c+u+\lambda v-ku_{0})],
\]
for any $\lambda\in(-\varepsilon,\varepsilon)$ and $k\geq0$. Since
$\mathcal{K}_{c}^{\varphi}$ is open, there exist a sufficiently small
$\alpha_{0}>0$ such that $u+\lambda v+\alpha|v|$ is in $\mathcal{K}_{c}^{\varphi}$
for all $\alpha\in[-\alpha_{0},\alpha_{0}]$. We can write
\[
\frac{g(\lambda+\alpha,k)-g(\lambda,k)}{\alpha}=\mathbb{E}\Bigl[\frac{1}{\alpha}\{\boldsymbol{\varphi}(c+u+(\lambda+\alpha)v-ku_{0})-\boldsymbol{\varphi}(c+u+\lambda v-ku_{0})\}\Bigr].
\]
The function in the expectation above is dominated by the $\mu$-integrable
function $\frac{1}{\alpha_{0}}\{\boldsymbol{\varphi}(c+u+\lambda v+\alpha_{0}|v|-ku_{0})-\boldsymbol{\varphi}(c+u+\lambda v-ku_{0})\}$.
By the Dominated Convergence Theorem, 
\begin{multline*}
\mathbb{E}\Bigl[\frac{1}{\alpha}\{\boldsymbol{\varphi}(c+u+(\lambda+\alpha)v-ku_{0})-\boldsymbol{\varphi}(c+u+\lambda v-ku_{0})\}\Bigr]\\
\rightarrow\mathbb{E}[v\boldsymbol{\varphi}'(c+u+\lambda v-ku_{0})],\qquad\text{as }\alpha\rightarrow0,
\end{multline*}
and, consequently,
\[
\frac{\partial g}{\partial\lambda}(\lambda,k)=\mathbb{E}[v\boldsymbol{\varphi}'(c+u+\lambda v-ku_{0})].
\]
Since $v\boldsymbol{\varphi}'(c+u+\lambda v-ku_{0})$ is dominated
by the $\mu$-integrable function $|v|\boldsymbol{\varphi}'(c+u+\varepsilon|v|-ku_{0})$,
we obtain for any sequence $\lambda_{n}\rightarrow\lambda$, 
\[
\mathbb{E}[v\boldsymbol{\varphi}'(c+u+\lambda_{n}v-ku_{0})]\rightarrow\mathbb{E}[v\boldsymbol{\varphi}'(c+u+\lambda v-ku_{0})],\quad\text{as }n\rightarrow\infty.
\]
Thus $\frac{\partial g}{\partial\lambda}(\lambda,k)$ is continuous
with respect to $\lambda$. Analogously, it can be shown that
\[
\frac{\partial g}{\partial k}(\lambda,k)=-\mathbb{E}[u_{0}\boldsymbol{\varphi}'(c+u+\lambda v-ku_{0})],
\]
and $\frac{\partial g}{\partial k}(\lambda,k)$ is continuous with
respect to $k$. The equality $g(\lambda,k(\lambda))=\mathbb{E}[\boldsymbol{\varphi}(c+u+\lambda v-k(\lambda)u_{0})]=1$
defines $k(\lambda)=\psi(u+\lambda v)$ as an implicit function of
$\lambda$. Notice that $\frac{\partial g(0,k)}{\partial k}<0$. By
the Implicit Function Theorem, the function $k(\lambda)=\psi(u+\lambda v)$
is continuously differentiable in a neighborhood of $0$, and has
derivative
\[
\frac{\partial k}{\partial\lambda}(0)=-\frac{(\partial g/\partial\lambda)(0,k(0))}{(\partial g/\partial k)(0,k(0))}.
\]
Consequently, 
\[
\partial\psi(u)(v)=\frac{\partial\psi(u+\lambda v)}{\partial\lambda}(0)=\frac{\mathbb{E}[v\boldsymbol{\varphi}'(c+u-\psi(u)u_{0})]}{\mathbb{E}[u_{0}\boldsymbol{\varphi}'(c+u-\psi(u)u_{0})]}.
\]
Thus the expression in \prettyref{eq:partial_psi} is the Gâteaux-derivative
of $\psi$.
\end{proof}

\begin{lem}
\label{lem:param_indep} Assume that $\varphi(t,\cdot)$ is continuously
differentiable. Then the divergence $D_{\psi}$ does not depend on
the parametrization of $\mathcal{F}_{c}^{\varphi}$.\end{lem}
\begin{proof}
For any $w\in\mathcal{B}_{c}^{\varphi}$, we denote $\widetilde{c}=c+w-\psi(w)u_{0}$.
Given $u,v\in\mathcal{B}_{c}^{\varphi}$, select $\widetilde{u},\widetilde{v}\in\mathcal{B}_{\widetilde{c}}^{\varphi}$
such that $\boldsymbol{\varphi}_{\widetilde{c}}(\widetilde{u})=\boldsymbol{\varphi}_{c}(u)$
and $\boldsymbol{\varphi}_{\widetilde{c}}(\widetilde{v})=\boldsymbol{\varphi}_{c}(v)$.
Let $\widetilde{\psi}\colon\mathcal{B}_{\widetilde{c}}^{\varphi}\rightarrow[0,\infty)$
be the normalizing function associated to $\widetilde{c}$. These
definitions provide 
\[
\widetilde{c}+\widetilde{u}-\widetilde{\psi}(\widetilde{u})u_{0}=c+u-\psi(u)u_{0},
\]
and 
\[
\widetilde{c}+\widetilde{v}-\widetilde{\psi}(\widetilde{v})u_{0}=c+v-\psi(v)u_{0}.
\]
Subtracting these equations, we obtain
\[
[-\widetilde{\psi}(\widetilde{v})+\widetilde{\psi}(\widetilde{u})]u_{0}+(\widetilde{v}-\widetilde{u})=[-\psi(v)+\psi(u)]u_{0}+(v-u)
\]
and, consequently, 
\begin{multline*}
\widetilde{\psi}(\widetilde{v})-\widetilde{\psi}(\widetilde{u})-\frac{\mathbb{E}[(\widetilde{v}-\widetilde{u})\boldsymbol{\varphi}'(\widetilde{c}+\widetilde{u}-\widetilde{\psi}(\widetilde{u})u_{0})]}{\mathbb{E}[u_{0}\boldsymbol{\varphi}'(\widetilde{c}+\widetilde{u}-\widetilde{\psi}(\widetilde{u})u_{0})]}\\
=\psi(v)-\psi(u)-\frac{\mathbb{E}[(v-u)\boldsymbol{\varphi}'(c+u-\psi(u)u_{0})]}{\mathbb{E}[u_{0}\boldsymbol{\varphi}'(c+u-\psi(u)u_{0})]}.
\end{multline*}
Therefore, $D_{\widetilde{\psi}}(\widetilde{u},\widetilde{v})=D_{\psi}(u,v)$.
\end{proof}

Let $p=\boldsymbol{\varphi}_{c}(u)$ and $q=\boldsymbol{\varphi}_{c}(v)$,
for $u,v\in\mathcal{B}_{c}^{\varphi}$. We denote the divergence between
the probability densities $p$ and $q$ by 
\[
D(p\mathbin{\Vert}q)=D_{\psi}(u,v).
\]
According to \prettyref{lem:param_indep}, $D(p\mathbin{\Vert}q)$
is well-defined if $p$ and $q$ are in the same $\varphi$-family.
We will find an expression for $D(p\mathbin{\Vert}q)$ where $p$
and $q$ are given explicitly. For $u=0$, we have $D(p\mathbin{\Vert}q)=D_{\psi}(0,v)=\psi(v)$,
and then
\[
D(p\mathbin{\Vert}q)=\frac{\mathbb{E}[(-v+\psi(v)u_{0})\boldsymbol{\varphi}'(c)]}{\mathbb{E}[u_{0}\boldsymbol{\varphi}'(c)]}.
\]
Therefore, the divergence between probability densities $p$ and $q$
in the same $\varphi$-family can be expressed as
\begin{equation}
D(p\mathbin{\Vert}q)=\frac{\mathbb{E}\biggl[\dfrac{\boldsymbol{\varphi}^{-1}(p)-\boldsymbol{\varphi}^{-1}(q)}{(\boldsymbol{\varphi}^{-1})'(p)}\biggr]}{\mathbb{E}\biggl[\dfrac{u_{0}}{(\boldsymbol{\varphi}^{-1})'(p)}\biggr]}.\label{eq:D_phi}
\end{equation}
Clearly, the expectation in \prettyref{eq:D_phi} may not be defined
if $p$ and $q$ are not in the same $\varphi$-family. We extend
the divergence in \prettyref{eq:D_phi} by setting $D(p\mathbin{\Vert}q)=\infty$
if $p$ and $q$ are not in the same $\varphi$-family. With this
extension, the divergence is denoted by $D_{\varphi}$ and is called
the \textit{$\varphi$-divergence}. By the strict convexity of $\varphi(t,\cdot)$,
we have the inequality $\varphi^{-1}(t,u)-\varphi^{-1}(t,v)\geq(\varphi^{-1})'(t,u)(u-v)$
for any $u,v>0$, with equality if and only if $u=v$. Hence $D_{\varphi}$
is always non-negative, and $D_{\varphi}(p\mathbin{\Vert}q)$ is equal
to zero if and only if $p=q$.

\begin{example}
With the variable $\kappa$-exponential $\exp_{\kappa}(t,u)=\exp_{\kappa(t)}(u)$
in the place of $\varphi(t,u)$, whose inverse $\varphi^{-1}(t,u)$
is the variable $\kappa$-logarithm $\ln_{\kappa}(t,u)=\ln_{\kappa(t)}(u)$,
we rewrite \prettyref{eq:D_phi} as
\begin{equation}
D(p\mathbin{\Vert}q)=\frac{\mathbb{E}\biggl[\dfrac{\operatorname{\mathbf{ln}}_{\boldsymbol{\kappa}}(p)-\operatorname{\mathbf{ln}}_{\boldsymbol{\kappa}}(q)}{\operatorname{\mathbf{ln}}_{\boldsymbol{\kappa}}'(p)}\biggr]}{\mathbb{E}\biggl[\dfrac{u_{0}}{\operatorname{\mathbf{ln}}_{\boldsymbol{\kappa}}'(p)}\biggr]},\label{eq:D_kappa}
\end{equation}
where $\operatorname{\mathbf{ln}}_{\boldsymbol{\kappa}}(p)$ denotes
$\ln_{\kappa(t)}(p(t))$. Since the $\kappa$-logarithm $\ln_{\kappa}(u)=\frac{u^{\kappa}-u^{-\kappa}}{2\kappa}$
has derivative $\ln_{\kappa}'(u)=\frac{1}{u}\frac{u^{\kappa}+u^{-\kappa}}{2}$,
the numerator and denominator in \prettyref{eq:D_kappa} result in
\[
\mathbb{E}\biggl[\dfrac{\operatorname{\mathbf{ln}}_{\boldsymbol{\kappa}}(p)-\operatorname{\mathbf{ln}}_{\boldsymbol{\kappa}}(q)}{\operatorname{\mathbf{ln}}_{\boldsymbol{\kappa}}'(p)}\biggr]=\mathbb{E}\Biggl[\dfrac{\dfrac{p^{\kappa}-p^{-\kappa}}{2\kappa}-\dfrac{q^{\kappa}-q^{-\kappa}}{2\kappa}}{\dfrac{1}{p}\dfrac{p^{\kappa}+p^{-\kappa}}{2}}\Biggr]=\frac{1}{\kappa}\mathbb{E}_{p}\Bigl[\dfrac{p^{\kappa}-p^{-\kappa}}{p^{\kappa}+p^{-\kappa}}-\dfrac{q^{\kappa}-q^{-\kappa}}{p^{\kappa}+p^{-\kappa}}\Bigr]
\]
and
\[
\mathbb{E}\biggl[\dfrac{u_{0}}{\operatorname{\mathbf{ln}}_{\boldsymbol{\kappa}}'(p)}\biggr]=\mathbb{E}_{p}\Bigl[\dfrac{2u_{0}}{p^{\kappa}+p^{-\kappa}}\Bigr],
\]
respectively. Thus \prettyref{eq:D_kappa} can be rewritten as
\[
D_{\kappa}(p\mathbin{\Vert}q)=\frac{1}{\kappa}\frac{\mathbb{E}_{p}\Bigl[\dfrac{p^{\kappa}-p^{-\kappa}}{p^{\kappa}+p^{-\kappa}}-\dfrac{q^{\kappa}-q^{-\kappa}}{p^{\kappa}+p^{-\kappa}}\Bigr]}{\mathbb{E}_{p}\Bigl[\dfrac{2u_{0}}{p^{\kappa}+p^{-\kappa}}\Bigr]},
\]
which we called the \textit{$\kappa$-divergence}.
\end{example}

\section*{Acknowledgments}

We are indebted to the anonymous referees for marked comments and
suggestions leading to the current version. This work received financial
support from CAPES -- Coordena\c{c}\~{a}o de Aperfei\c{c}oamento
de Pessoal de N\'{\i}vel Superior.

\appendix
\bibliographystyle{plain}
\bibliography{refs_phi-families}

\end{document}